\documentclass{amsart}
\usepackage{amsmath,amssymb}

\usepackage{epsfig}
\usepackage{url}
\usepackage{amssymb}
\usepackage{amsmath}
\usepackage{amsfonts}

\def\Dbar{\leavevmode\lower.6ex\hbox to 0pt{\hskip-.23ex \accent"16\hss}D}

\def\bZ{{\mbox{\bf Z}}}

\def\bX{{\mbox{\bf X}}}

\begin{document}

\title{Two classes of Hadamard matrices of Goethals-Seidel type}
\author {Dragomir {\v{Z}}. {\Dbar}okovi{\'c}}
\address{University of Waterloo, 
Department of Pure Mathematics,
Waterloo, Ontario, N2L 3G1, Canada}
\email{dragomir@rogers.com}
\date{}

\begin{abstract}
We introduce two classes of Hadamard matrices of
Goethals-Seidel type and construct many matrices in these classes.
The largest among them have order $4\cdot 631=2524$.
We do not know whether these classes are infinite. 

No skew Hadamard matrix of order $4\cdot 397$ was known so far.
We constructed eight such matrices, two of them are 
listed in Section \ref{Additional}.
\end{abstract}

\maketitle

\section{Introduction}

The Goethals-Seidel (GS) array 
\begin{equation*} 
\left[ \begin{array}{cccc}
A_0    & A_1 R & A_2 R & A_3 R \\
-A_1 R & A_0   & -RA_3 & RA_2  \\
-A_2 R & RA_3  & A_0   & -RA_1 \\
-A_3 R & -RA_2 & RA_1  & A_0  \\
\end{array} \right]
\end{equation*} 
is a powerful tool for constructing Hadamard matrices 
\cite{DK:Goethals-Seidel,GS-array:1970,SY:1992} and in particular
skew-Hadamard matrices.

In order to use this tool we need a special kind of 
difference families over a finite abelian group $G$ 
consisting of four blocks
\begin{equation} \label{dfm}
\bX=(X_0,X_1,X_2,X_3).
\end{equation}
Let the parameter set of this family be
\begin{equation} \label{prm}
(v; k_0,k_1,k_2,k_3; \lambda).
\end{equation}
Thus we have $v=|G|$, $k_i=|X_i|$ for $0\le i\le 3$, and
\begin{equation} \label{prva}
	\sum_{i=0}^3 k_i(k_i-1)=\lambda(v-1).
\end{equation}
An additional parameter, called {\em order}, is defined by
\begin{equation} \label{druga}
	n=\sum_{i=0}^3 k_i-\lambda.
\end{equation}
If $n=v$ we say that the difference family $\bX$ and 
its parameter set are of {\em Goethals-Seidel (GS) type}.

Given such a family, there is a simple procedure which uses the
GS-array to produce a Hadamard matrix of order $4v$. We say that
such Hadamard matrix has {\em GS-type}. If the group $G$ is
cyclic then the matrix blocks 
$A_i$ are circulants of order $v$ and $R$ is the back-circulant
identiy matrix, for more details see 
\cite{DK:Goethals-Seidel,SY:1992}.

It is an open question whether for each GS-parameter set 
\eqref{prm} there exists a GS-difference family
(over some abelian group $G$ of order $v$) having these parameters.
A positive answer would imply that the well known
Hadamard matrix conjecture is true. It is a folklore conjecture
that the answer to the above question is positive. There is also
a stronger conjecture which requires $G$ to be a cyclic group.

Assume now that $\bX$ is a GS-difference family
and let $i\in\{0,1,2,3\}$.
If we replace the block $X_i$ with its
complement $G \setminus X_i$ we get again a GS-difference 
family but now its parameter set is changed: we have to replace
$k_i$ with $v-k_i$ and replace $\lambda$ with $\lambda+v-2k_i$. 
Thus, in order to avoid redundancy,
we may and we shall assume that each $k_i \le v/2$.

Since the size of a Hadamard matrix can easily be doubled, 
we consider only the cases with $v$ odd. Moreover, we shall assume
that $G$ is cyclic and we identify $G$ with the additive group of
the ring $\bZ_v$ (integers modulo $v$). 
By $\bZ_v^*$ we denote the group of units of the ring $\bZ_v$.

We say that a GS-parameter set \eqref{prm} is {\em special} if $k_1=k_2=k_3$.
By using the equations \eqref{prva},  \eqref{druga} and $n=v$, 
one can easily show that if a parameter set $(v;k_0,k,k,k;\lambda)$ is
special then the integer $1+2(k_0+3k)-3(k_0-k)^2$ must be a square.

Our objective is to construct pairs $(\bX,\mu)$ where $X$ is a 
GS-difference family and $\mu\in\bZ_v^*$ (a {\em multiplier}) 
sends $X_1$ to $X_2$ and $X_2$ to $X_3$, i.e.,
$\mu X_1=X_2$ and $\mu X_2=X_3$. There are two cases:
$\mu X_3=X_1$ and $\mu X_3\ne X_1$. In the first case we say that 
$(\bX,\mu)$ and the associated Hadamard matrix have a {\em spin structure}
and in the second case a {\em slide structure}.

Note that if $(\bX,\mu)$ is such a pair then so is the pair
$(\bX',\mu^{-1})$ where $\bX'$ is obtained from $\bX$ by interchanging the 
blocks $X_1$ and $X_3$.

Note also that if $(\bX,\mu)$ has spin structure then 
$\bZ_v^*$ must have an element of order 3.
For convenience we say that a special parameter set is 
of {\em spin type} if $\bZ_v^*$ has an element of order 3.

We point out that there exist pairs
$(\bX,\mu_1)$ and $(\bX,\mu_2)$ with $\mu_1 \ne \mu_2$ 
such that $(\bX,\mu_1)$ has spin structure while 
$(\bX,\mu_2)$ has slide structure. 
We shall give an example at the end of Section \ref{Additional}.

A subset $X$ of $\bZ_v$ is {\em symmetric} if $-X=X$. It is 
{\em skew} if $\bZ_v$ is a disjoint union of 
$X$, $-X$ and $\{0\}$. If the block $X_0$ of a GS-difference family $\bf X$ is skew then
the associated Hadamard matrix is of skew type.
We say that a GS-difference family $\bX$ is {\em good} if $X_0$ is skew
and the other three blocks are symmetric.
On the other hand, if $X_0$ is symmetric and the other three blocks are skew
then we say that $\bX$ is a {\em best} family.
(The corresponding matrix blocks $A_i$ in the GS-array have the same symmetry 
type as the corresponding $X_i$.)

From now on we assume that $\bX$ is a GS-difference family whose parameter set is special. 
Hence the blocks $X_1,X_2,X_3$ have the same size and the same symmetry type. 

We indicate the symmetry types of $X_0$ and $X_1$ by letters,
$s$ for symmetric and $k$ for skew. For instance
the symmetry symbol $(sk)$ means that $X_0$ is symmetric and 
$X_1$ is skew (and so are $X_2$ and $X_3$).
The meaning of $(s\ast)$ is that $X_0$ is symmetric and that
we make no claim about the symmetry type of $X_1$.
The symmetry symbol $(\ast\ast)$ is usually omitted.
GS-type Hadamard matrices having special parameter sets and symmetry type $(ss)$, $(sk)$, $(ks)$
are special cases of Williamson, good and best matrices, respectively.
For the list of known Williamson matrices with $v$ odd and $v\le59$ see \cite{HKT}.

\section {Parameters and Hadamard matrices of spin type for $v<100$} \label{spin}

In this section we list all parameter sets of spin type with $v<100$. 
We give examples of spin difference families whenever we succeeded to find one.
In some cases we found many solutions but we list only a few.
The question mark after the parameter set means that we failed to find an associated
spin difference family. In such cases we searched for the slide difference families
(and if found, they are given in Section \ref{slide}).
Beside the blocks $X_0$ and $X_1$ we also specify the multiplier $\mu$.
In general $\mu$ does not have to fix the block $X_0$.
If it does, we point it out by saying that $X_0$ is fixed.

A few examples in our list below are extracted from other papers or some well known lists of
special Hadamard matrices such as Williamson matrices, good matrices and best matrices. 
We point out explicitly such cases. \vspace{2mm}

$ \begin{array}{ll}
	(7; 3,2,2,2; 2)    & (ks),~ \mu=2,~ X_0~ {\rm~ fixed}, \\
	X_0=[1,2,4],~ X_1=[1,6] & \\
	{\rm Good~ matrices.}   &
\end{array} $ \vspace{2mm}

$ \begin{array}{ll}
	(7; 1,3,3,3; 3)      & (ss),~ \mu=2,~ X_0 {\rm~ fixed}, \\
	X_0=[0],~ X_1=[0,1,6]       & \\
	{\rm Williamson~ matrices.} &
\end{array} $ \vspace{2mm}

$ \begin{array}{ll}
	(9; 3,3,3,3; 3)     & (ss),~ \mu=4,~ X_0~ {\rm fixed}, \\
	X_0=[0,3,6],~ X_1=[0,1,8]   & \\
	{\rm Williamson~ matrices.} &
\end{array} $ \vspace{2mm}

$ \begin{array}{ll}
	(13; 4,5,5,5; 6) & (\ast s),~ \mu=3, \\
	X_0=[0,1,4,6]    &  X_1=[0,4,6,7,9]. 
\end{array} $ \vspace{2mm}

$ \begin{array}{ll}
	(13; 3,6,6,6; 8) & (\ast k),~ \mu=3, \\
	X_0=[6,7,10]     & X_1=[1,2,3,4,6,8].
\end{array} $ \vspace{2mm}

$ \begin{array}{ll}
	(19; 9,7,7,7; 11)    & (ks),~ \mu=-2,~ X_0~ {\rm fixed}, \\
	X_0=[1,4,5,6,7,9,11,16,17] & X_1=[0,1,7,8,11,12,18]. \\
	{\rm Good~ matrices.} &
\end{array} $ \vspace{2mm}

In the next two examples the block $X_0$ can be combined
with each of the listed choices for $X_1$. \vspace{2mm}

$ \begin{array}{ll}
	(19; 9,7,7,7; 11) & (k\ast),~ \mu=7,~ X_0~ {\rm fixed},\\
	X_0=[2,3,4,6,8,9,12,14,18] & \\
	X_1=[0,1,6,8,13,15,16],    & X_1=[0,2,8,9,10,11,13], \\
	X_1=[0,2,8,9,10,11,14],    & X_1=[0,1,4,5,12,15,18]. 
\end{array} $    \vspace{2mm}

$ \begin{array}{ll}
	(19; 6,8,8,8; 11), & (s\ast),~ \mu=7,~ X_0~ {\rm fixed}, \\
	X_0=[4,6,9,10,13,15]        & X_1=[0,1,5,8,9,10,11,13], \\
	X_1=[0,1,8,9,10,11,13,14],  & X_1=[0,1,8,9,10,11,14,17], \\
	X_1=[0,4,7,9,10,11,14,18],  & X_1=[0,4,8,9,10,11,12,13], \\
	X_1=[0,6,9,10,11,12,14,18], & X_1=[0,6,9,10,11,12,16,17].
\end{array} $ \vspace{2mm}

$ \begin{array}{ll}
	(21; 9,8,8,8; 12)            & (\ast s),~ \mu=4, \\
	X_0=[1,4,5,8,10,11,12,17,19] & X_1=[1,3,8,9,12,13,18,20].
\end{array} $ \vspace{2mm}

$ \begin{array}{ll}
	(21; 6,10,10,10; 15)  & \mu=4,\\
	X_0=[0,8,10,11,12,16] & X_1=[0,1,2,5,6,7,13,15,16,19].
\end{array} $ \vspace{2mm}

Although the next parameter set is of spin type,
our exhaustive search did not find any (cyclic) spin
difference family having these parameters. \vspace{2mm}

$ \begin{array}{c}
(27; 9,12,12,12; 18).
\end{array} $ \vspace{2mm}

$ \begin{array}{l}
	(31; 12,13,13,13; 20),~ (\ast s),~ \mu=5,\\
X_0=[2,4,6,12,14,16,17,19,25,26,28,29],\\
X_1=[0,3,9,11,13,14,15,16,17,18,20,22,28].
\end{array} $ \vspace{2mm}

$ \begin{array}{l}
(31; 10,15,15,15; 24),~ (\ast k),~ \mu=5, \\
X_0=[0,2,10,13,16,17,20,26,28,29], \\
X_1=[1,4,12,13,14,15,20,21,22,23,24,25,26,28,29].
\end{array} $ \vspace{2mm}

$ \begin{array}{l}
(37; 18,15,15,15; 26),~ (ks),~ \mu=10,~ X_0~ {\rm fixed}, \\
X_0=[2,3,4,6,8,11,15,18,20,21,23,24,25,27,28,30,32,36],\\
X_1=[0,1,2,5,9,13,14,15,22,23,24,28,32,35,36].\\
{\rm Good~ matrices.}
\end{array} $ \vspace{2mm}

$ \begin{array}{l}
(37; 13,17,17,17; 27),~ (ss),~ \mu=10,~ X_0~ {\rm fixed}, \\
X_0=[0,3,4,5,7,13,18,19,24,30,32,33,34], \\
X_1=[0,1,2,3,4,6,12,13,18,19,24,25,31,33,34,35,36]. \\
{\rm Williamson~ matrices.} 
\end{array} $ \vspace{2mm}

$ \begin{array}{l}
(39; 18,16,16,16; 27),~ (s\ast),~ \mu=16,~ X_0~ {\rm fixed},\\
X_0=[1,4,6,10,14,15,16,17,18,21,22,23,24,25,29,33,35,38],\\
X_1=[0,2,4,6,7,10,11,14,16,19,20,22,26,32,33,38].
\end{array} $ \vspace{2mm}

$ \begin{array}{l}
(39; 15,17,17,17; 27),~ (ss),~ \mu=16,~ X_0~ {\rm fixed},\\
X_0=[0,4,8,10,11,13,14,19,20,25,26,28,29,31,35],\\
X_1=[0,1,4,5,6,8,11,12,14,25,27,28,31,33,34,35,38].\\
{\rm Williamson~ matrices.}
\end{array} $ \vspace{2mm}

$ \begin{array}{l}
(43; 19,18,18,18; 30)~ ?\\
(43; 15,21,21,21; 35)~ ?\\
(49; 21,21,21,21; 35)~ ?\\
(49; 19,22,22,22; 36)~ ?\\
(57; 24,25,25,25; 42)~ ?
\end{array} $ \vspace{2mm}

An exhaustive search for best matrices of order 57 has been
carried out in \cite{BDKG:best2019}. Three solutions were found.
It turned out that two of them have spin structure: \vspace{2mm}

$ \begin{array}{l}
(57; 21,28,28,28; 48),~ (sk),~ \mu=7,~ X_0~ {\rm is~ fixed~ in~ both~ cases.}\\
1)~ X_0=[0,4,11,12,18,19,20,25,26,27,28,29,30,31,32,37,38,39,45,46,53],\\
\hspace{4mm}   X_1=[1,2,5,6,7,8,9,10,12,17,19,21,22,24,25,28,30,31,34,\\
\hspace{13mm}           37,39,41,42,43,44,46,53,54],\\
2)~ X_0=[0,4,5,12,17,18,19,22,25,27,28,29,30,32,35,38,39,40,45,52,53],\\
\hspace{4mm}   X_1=[2,4,6,7,8,10,13,16,17,18,19,20,22,23,24,25,26,27,28, \\
\hspace{13mm}           36,42,43,45,46,48,52,54,56].\\
{\rm Best~ matrices.}
\end{array} $ \vspace{2mm}

When all blocks $X_i$ are $H$-invariant for some 
subgroup $H$ of $\bZ_v^*$ then it suffices to list only 
the representatives of the $H$-orbits contained in 
$X_0$ (rep. 0) and $X_1$ (rep. 1). \vspace{2mm} 

$ \begin{array}{l}
(61; 30,26,26,26; 47),~ (k\ast),~ \mu=13,~ H=[1,9,20,34,58], \\
   {\rm rep.}~ 0:~ [3,4,5,6,8,10],~ {\rm rep.}~ 1:~ [0,8,10,13,23,26].\\
   X_0~ {\rm fixed.}
\end{array} $ \vspace{2mm}

$ \begin{array}{l}
(61; 24,28,28,28; 47)~ ? \\
(63; 30,27,27,27; 48)~ ? \\
(63; 24,30,30,30; 51)~ ? \\
(67; 31,29,29,29; 51)~ ? \\
(67; 28,30,30,30; 51)~ ? \\
(73; 33,32,32,32; 56)~ ?
\end{array} $ \vspace{2mm}

$ \begin{array}{l}
(73; 28,36,36,36; 63),~ (\ast s),~ \mu=4,~ H=[1,8,64],\\
1)~ {\rm rep.}~ 0:~ [0,9,13,18,25,26,27,35,36,43],\\
	\hspace{4mm}   {\rm rep.}~ 1:~ [1,2,4,9,11,14,18,21,26,34,36,43].\\
2)~ {\rm rep.}~ 0:~ [0,5,7,9,14,17,18,33,34,36],\\
	\hspace{4mm}   {\rm rep.}~ 1:~ [4,5,7,12,14,17,21,33,34,35,36,43].\\
X_0~ {\rm is~ fixed~ in~ both~ cases.}
\end{array} $ \vspace{2mm}

$ \begin{array}{l}
(79; 33,36,36,36; 62),~ \mu=12,~ H:=[1,23,55], \\
1)~ {\rm rep.}~ 0:~ [[2,9,12,15,20,22,24,34,37,40,47],~ (\ast s),\\
        \hspace{4mm} {\rm rep.}~ 1:~ [4,5,6,12,17,18,20,22,37,40,41,47].\\
2)~ {\rm rep.}~ 0:~ [4,5,6,10,17,22,27,30,33,44,47],\\
        \hspace{4mm} {\rm rep.}~ 1:~ [1,5,8,17,18,20,22,34,37,40,44,47].\\
\end{array} $ \vspace{2mm}

$ \begin{array}{l}
(81; 36,36,36,36; 63)~ ? 
\end{array} $ \vspace{2mm}

$ \begin{array}{l}
(91; 45,40,40,40; 74),~ \mu=9,~ H=[1,16,74],\\
   {\rm rep.}~ 0: [1,5,6,7,13,14,16,19,20,24,28,29,39,47,49],\\
   {\rm rep.}~ 1: [0,3,8,13,16,23,24,38,40,46,47,48,49,57].
\end{array} $ \vspace{2mm}

The three difference families listed below are pairwise inequivalent although 
the last two share the same block $X_0$.

$ \begin{array}{l}
(91; 40,41,41,41; 72),~ \mu=9,~ H=[1,22,29], \\
1)~ {\rm rep.}~ 0:~ [0,3,4,5,6,10,12,23,27,31,36,39,40,62,65,78],\\
   \hspace{4mm} {\rm rep.}~ 1:~ [2,4,5,8,11,20,23,24,27,28,36,39,40,49,52].\\
2)~ {\rm rep.}~ 0:~ [0,1,3,4,9,20,23,24,27,33,36,39,40,53,65,78],\\
   \hspace{4mm} {\rm rep.}~ 1:~ [6,9,14,15,18,23,24,27,28,33,36,37,39,40,52].\\
3)~ {\rm rep.}~ 0:~ [0,1,3,4,9,20,23,24,27,33,36,39,40,53,65,78],\\
   \hspace{4mm} {\rm rep.}~ 1:~ [1,2,5,7,15,20,23,24,27,28,36,39,40,52,53].
\end{array} $ \vspace{2mm}

$ \begin{array}{l}
(91; 37,43,43,43; 75), \\
1)~ (s\ast),~ \mu=9,~ H=[1,16,74],\\
   {\rm rep.}~ 0:~ [0,3,4,5,8,11,19,25,27,43,45,50,55],\\
   {\rm rep.}~ 1:~ [0,1,4,5,13,14,15,25,28,33,38,43,44,49,55].\\
2)~ \mu=16,~ H=[1,9,81],\\
   {\rm rep.}~ 0:~ [0,6,8,12,13,19,20,24,38,39,40,48,57],\\
   {\rm rep.}~ 1:~ [0,2,15,16,19,23,24,28,30,38,40,47,48,49,57].
\end{array} $ \vspace{2mm}

$ \begin{array}{l}
(91; 36,45,45,45; 80),~ \mu=12,~ H=[1,9,81], \\
   {\rm rep.}~ 0:~ [2,7,10,12,13,14,15,28,38,39,49,57],\\
   {\rm rep.}~ 1:~ [4,14,15,16,19,20,23,24,28,38,46,47,48,49,57].
\end{array} $ \vspace{2mm}

$ \begin{array}{l}
(93; 45,41,41,41; 75),~ \mu=25,~ H=[1,4,16,64,70], \\
   {\rm rep.}~ 0:~ [3,10,11,14,21,23,33,34,46],\\
   {\rm rep.}~ 1:~ [3,9,11,17,23,33,34,46,62].
\end{array} $ \vspace{2mm}

$ \begin{array}{l}
(93; 39,43,43,43; 75)~ ? \\
(97; 46,43,43,43; 78)~ ? \\
(97; 39,47,47,47; 83)~ ?
\end{array} $

\section{Spin-type GS-difference families with $v>100$} \label{Additional}

There are only finitely many known good matrices. Their orders
are $\le 127$. We constructed those of order 127 long ago in 
\cite{Djokovic:JCMCC:1993}. Much later we observed that the
GS-difference family used in that construction
has very special properties to which we now refer as
spin structure. Here is that example
where we replaced the three blocks of size 70 by their
complements (of size 57): \vspace{2mm}

$ \begin{array}{l}
(127; 63,57,57,57; 107),~ (ks),~ \mu=19,~ H=[1,2,4,8,16,32,64], \\
   {\rm rep.}~ 0:~ [1,3,7,9,11,19,21,23,47], \\
   {\rm rep.}~ 1:~ [0,3,7,9,11,15,29,31,55]. \\
	{\rm Good~ matrices}, X_0~ {\rm fixed.}
\end{array} $ \vspace{2mm}

The $H$ above is the subgroup of $\bZ_{127}^*$ of order 7.
All four blocks are $H$-invariant, and so each of them is a 
union of some $H$-orbits. The representatives of the $H$-orbits
contained in $X_0$ and those in $X_1$ are listed above. 
The block $X_0$ is skew and $X_1$ is symmetric.
The multiplier $\mu=19$ (of order 3) fixes the block $X_0$
and permutes cyclically the other three blocks.

As $v$ grows the search for spin Hadamard matrices becomes
harder and harder. The search may be feasible only for difference families
invariant under relatively large subgroup
say $H$ of $\bZ_v^*$. In the following examples we use
subgroups of order 7,9,11,13,15 and 27. \vspace{2mm}

$ \begin{array}{l}
(129; 63,58,58,58; 108),~ \mu=13,~ H=[1,4,16,64,97,121,127],\\
   {\rm rep.}~ 0:~ [1,9,10,14,19,21,23,26,27], \\
   {\rm rep.}~ 1:~ [2,5,9,10,13,18,22,27,43,86].
\end{array} $ \vspace{2mm}

$ \begin{array}{l}
(271; 135,126,126,126; 242),~ (k\ast),~ \mu=5, \\
H=[1,28,106,125,169,178,242,248,258], \\
1)~ {\rm rep.}~ 0:~ [1,4,5,7,8,11,14,16,19,21,22,25,31,43,44],\\
	\hspace{4mm}   {\rm rep.}~ 1:~ [1,2,3,5,7,8,12,19,22,27,38,42,44,51],\\
2)~ {\rm rep.}~ 0:~ [1,2,4,5,7,8,9,11,14,16,17,22,25,31,44],\\
	\hspace{4mm} {\rm rep.}~ 1:~ [1,3,5,9,12,14,17,19,21,22,33,44,71,86].
\end{array} $ \vspace{2mm}

$ \begin{array}{l}
(331; 165,155,155,155; 299),~ (k\ast),~ \mu=31,\\
H=[1,74,80,85,111,120,167,180,270,274,293], \\
1)~ {\rm rep.}~ 0:~ [5,10,11,13,16,19,20,22,32,38,53,56,64,76,101],\\
	\hspace{4mm}    {\rm rep.}~ 1:~ [0,4,11,16,20,28,31,37,41,49,53,56,73,88,101].\\
2)~ {\rm rep.}~ 0:~ [4,5,13,14,16,19,20,22,32,38,49,53,56,64,76],\\
	\hspace{4mm}    {\rm rep.}~ 1:~ [0,11,13,14,19,22,31,37,44,49,56,62,73,76,88].
\end{array} $ \vspace{2mm}

No skew Hadamard matrix of order $4\cdot 397$
is known, see e.g. \cite[Table 9.2]{SY:2020},
 \cite[Table 1]{KouStyl:2008} or \cite[Table 1.51]{CK-Had:2007}.
 We constructed eight such matrices from spin-type GS-difference families of symmetry type $(k\ast)$.
 Two of them are presented below.  \vspace{2mm}

$ \begin{array}{l}
(397; 198,187,187,187; 362),~ \mu=34, \\ 
H=[1,16,31,99,126,167,256,273,290,333,393], \\
1)~ {\rm rep.}~ 0:~ [1,6,7,8,9,10,11,12,17,18,20,21,29,34,46,47,53,106],~ (k\ast),\\
\hspace{4mm} {\rm rep.}~ 1:~ 
[2,11,12,17,18,20,24,27,33,34,36,40,46,47,53,58,71].\\
2)~ {\rm rep.}~ 0:~ [3,4,6,7,8,17,18,20,21,23,27,30,33,36,40,44,46,53],~ (k\ast),\\
\hspace{4mm} {\rm rep.}~ 1:~ 
[1,8,12,17,18,20,24,27,33,34,40,44,46,47,53,58,71]. \\
3)~ {\rm rep.}~ 0:~ [3,5,9,10,11,12,18,20,21,23,29,33,36,40,44,47,61,72],
~ (s\ast), \\
\hspace{4mm} {\rm rep.}~ 1:~ [2,3,6,10,17,22,24,33,34,36,40,46,47,53,58,71,72]. \\
4)~ {\rm rep.}~ 0:~ [1,2,4,6,7,8,15,17,22,24,27,30,34,46,53,58,71,106],~ (s\ast),\\
\hspace{4mm} {\rm rep.}~ 1:~ [3,5,6,11,15,17,20,22,23,27,33,34,46,47,53,58,71].\\
X_0~ {\rm fixed~ in~ all~ four~ cases.}
\end{array} $ \vspace{2mm}

For $v=547$ Paley's construction gives a skew Hadamard matrix of order $4v$. 
We computed eight spin-type skew Hadamard matrices of the same order.
The blocks $X_0$ and $X_1$ for one of them are specified below:  \vspace{2mm}

$ \begin{array}{l}
(547; 273,260,260,260; 506),~ (k\ast),~ \mu=40, \\
H=[1,46,237,261,293,350,353,375,440,475,509,517,519], \\
{\rm rep.}~ 0:~ [1,4,5,6,10,11,13,14,17,25,29,34,35,40,49,52,55,64,69,110,123],\\
{\rm rep.}~ 1:~ [1,4,5,11,16,17,20,26,32,33,34,41,49,52,55,64,70,80,123,207].\\
X_0~ {\rm fixed.}
\end{array} $ \vspace{2mm}

So far no spin-type GS-difference families with $v>631$ are known.
We list below five spin-type GS-difference families that we constructed for $v=631$. 
The first two are just minor modifications of the two GS-difference families used in \cite{DGK:JCD:2014}
to construct two skew-Hadamard matrices of order $4v$. 
Surprisingly, it turned out that these two families indeed have spin structure.
In all five families, $H$ is the subgroup of order 15 and the block $X_0$ is skew 
and fixed by $\mu$. \vspace{2mm}
 
$ \begin{array}{l}
(631; 315,301,301,301; 587),~ (k\ast),~ \mu=2, \\
H=[1,8,43,64,79,188,228,242,279,310,339,344,512,562,587], \\
	1)~ {\rm rep.}~ 0:~ [1,2,3,4,6,7,12,13,14,17,19,21,26,27,31,38,\\ \hspace{2cm}	42,52,62,76,124],\\
	\hspace{5mm}	{\rm rep.}~ 1:~ [0,11,13,14,18,19,21,22,29,35,39,46,62,63,\\ \hspace{2cm}	65,66,67,92,117,124,187],\\
	2)~ {\rm rep.}~ 0:~ [11,13,19,22,26,29,31,33,38,39,44,52,62,65,\\ \hspace{2cm}	66,67,76,78,117,124,187],\\
	\hspace{5mm}	{\rm rep.}~ 1:~ [0,2,6,7,12,13,19,21,27,31,35,44,52,63,66,\\ \hspace{2cm}	76,78,92,124,126,187],\\
	3)~ {\rm rep.}~ 0:~ [1,2,4,5,7,9,14,17,18,21,23,27,31,33,42,46,\\ \hspace{2cm}	62,66,67,92,124],\\
	\hspace{5mm}	{\rm rep.}~ 1:~ [0,4,5,9,13,14,19,21,22,27,29,31,33,35,44,\\ \hspace{2cm}	63,76,92,124,126,187],\\
	4)~ {\rm rep.}~ 0:~ [1,2,4,7,13,14,23,26,27,31,33,46,52,62,65,\\ \hspace{2cm}	66,67,92,117,124,187],\\
	\hspace{5mm}	{\rm rep.}~ 1:~ [0,3,4,6,12,13,17,18,19,23,26,27,29,31,35,\\ \hspace{2cm}	46,62,65,67,76,92],\\
	5)~ {\rm rep.}~ 0:~ [1,2,3,4,6,7,12,13,14,17,19,21,26,27,31,38,\\ \hspace{2cm}	42,52,62,76,124],\\
	\hspace{5mm}	{\rm rep.}~ 1:~ [0,2,5,6,7,13,18,19,21,27,33,39,44,52,62,\\ \hspace{2cm}	63,76,78,92,117,126].
\end{array} $ \vspace{2mm}

If we replace the multiplier $\mu=2$ with $\mu=4$, then all five pairs $(\bX,\mu)$
above will loose the spin structure and acquire the slide structure.

\section {Hadamard matrices with slide structure} \label{slide}

In this section we give examples of pairs $(\bX,\mu)$ having
slide structure. In three of these examples we have indicated that
$X_3=-X_1$. Hence we can replace $X_3$ with $X_1$ to obtain 
a GS-difference family with a repeated block. \vspace{2mm}

$ \begin{array}{l}
(25; 10,10,10,10; 15),~ \mu=7, \\
   X_0=[1,5,7,9,12,13,15,16,19,21],\\
   X_1=[2,4,5,6,7,11,19,20,21,24].\\
X_3 = -X_1.
\end{array} $ \vspace{2mm}

$ \begin{array}{l}
(27; 9,12,12,12; 18),~ \mu=7, \\
1)~ X_0=[2,6,8,10,13,14,16,18,21],\\
\hspace{4mm}   X_1=[2,3,8,11,12,13,14,15,16,17,22,26].\\
2)~ X_0=[2,6,8,9,13,14,16,19,21],\\
\hspace{4mm}   X_1=[1,5,7,10,12,13,14,15,16,17,21,25].
\end{array} $ \vspace{2mm}

$ \begin{array}{l}
(31; 10,15,15,15; 24),~ \mu=4,~ H=[1,5,25],\\
   {\rm rep.}~ 0:~ [0,3,11,12],~ {\rm rep.}~ 1:~ [1,3,8,16,17].
\end{array} $ \vspace{2mm}

$ \begin{array}{l}
(43; 19,18,18,18; 30),~ \mu=3,~ H=[1,6,36],\\
   {\rm rep.}~ 0:~ [0,1,3,9,19,20,21],~ {\rm rep.}~ 1;~ [1,5,13,19,20,26].
\end{array} $ \vspace{2mm}

$ \begin{array}{l}
(43; 15,21,21,21; 35),~ H=[1,6,36],\\
1)~ {\rm rep.}~ 0:~ [1,4,9,20,26],~ {\rm rep.}~ 1:~ [1,9,10,13,20,21,26],~ \mu=3,\\
2)~ {\rm rep.}~ 0:~ [1,2,3,20,21],~ {\rm rep.}~ 1:~ [3,4,5,10,20,21,26],~ \mu=9.
\end{array} $ \vspace{2mm}

$ \begin{array}{l}
(49; 21,21,21,21; 35),~ \mu=9,~ H=[1,18,30],\\
   {\rm rep.}~ 0:~ [1,2,8,9,21,24,29],~ {\rm rep.}~ 1:~ [2,6,12,19,24,26,29].
\end{array} $ \vspace{2mm}

$ \begin{array}{l}
(49; 19,22,22,22; 36),~ \mu=4,~ H=[1,18,30], \\
   {\rm rep.}~ 0:~ [0,2,4,7,8,16,29],~ {\rm rep.}~ 1:~ [0,2,7,8,12,13,26,29].
\end{array} $ \vspace{2mm}
 
$ \begin{array}{l}
(61; 30,26,26,26; 47),~ (s\ast),~ \mu=2,~ H=[1,9,20,34,58], \\
   {\rm rep.}~ 0:~ [1,3,4,5,12,13],~ {\rm rep.}~ 1:~ [0,8,12,13,23,26].
\end{array} $ \vspace{2mm}

$ \begin{array}{l}
(61; 24,28,28,28; 47),~ \mu=11,~ H=[1,13,47],\\
   {\rm rep.}~ 0:~ [1,8,9,11,12,18,27,36], \\  
   {\rm rep.}~ 1:~ [0,1,7,8,9,18,22,27,31,32].\\
   X_3 = -X_1.
\end{array} $ \vspace{2mm}

$ \begin{array}{l}
(67; 28,30,30,30; 51),~ \mu=15,~ H=[1,29,37],\\
   {\rm rep.}~ 0:~ [0,1,4,5,7,11,12,16,17,19],\\
   {\rm rep.}~ 1:~ [2,3,4,9,12,15,16,17,19,21].
\end{array} $ \vspace{2mm}

$ \begin{array}{l}
(73; 28,36,36,36; 63),~ \mu=3,~ H=<2>~ 
	{\rm subgroup~ of~ order~ 9},\\
   {\rm rep.}~ 0:~ [0,1,3,13],~ {\rm rep.}~ 1:~ [3,11,13,17],\\
   X_3 = -X_1.
\end{array} $ \vspace{2mm}

$ \begin{array}{l}
(129; 63,58,58,58; 108),~ \mu=19,~ H=[1,4,16,64,97,121,127],\\
   {\rm rep.}~ 0:~ [1,6,9,11,13,19,23,26,27], \\
   {\rm rep.}~ 1:~ [2,7,9,14,21,22,23,27,43,86].
\end{array} $ \vspace{2mm}

$ \begin{array}{l}
(211; 91,105,105,105; 195),~ H=<19> 
{\rm subgroup~ of~ order~ 15,}\\
1)~ {\rm rep.}~ 0:~ [0,1,2,10,22,26,43],~ \mu=25,\\
	\hspace{4mm}  {\rm rep.}~ 1:~ [1,4,10,11,22,29,43].\\
2)~ {\rm rep.}~ 0:~ [5,10,11,22,26,29],~ \mu=13,\\
	\hspace{4mm}   {\rm rep.}~ 1:~ [2,5,7,10,22,29,43].
\end{array} $ \vspace{2mm}

\section{Appendix: Prime chains}

While trying to construct Hadamard matrices of unknown small
orders $4v$, we noticed that the hardest cases are often those where $v$ is a prime congruent to 3 mod 4 and $(v-1)/2$ is also 
a prime. Let us mention three of them. For about 30 years (1932-1962) it was not known how to construct a Hadamard matrix 
of order $92=4\cdot 23$, 
see \cite[p. 177]{LW:Combinatorics:1992}. 
It was the smallest order (multiple of 4) for which no Hadamard matrix was known.
Such matrix was constructed in 1962 when the computers became available.
For about 20 years (1985-2005) the smallest unknown order was $428=4\cdot 107$.
The first Hadamard matrix of that order was constructed by H. Kharaghani and B. Tayfeh-Rezaie \cite{KT:428:2005}.
Presently, for $v<250$ there are only three values, namely 
$v=167,179,223$ for which no Hadamard matrix of order $4v$ is known.
All three are primes congruent to 3 mod 4 and in the first two cases the number $(v-1)/2$ is a prime.
These observations led us to consider some questions about the prime numbers.

Let $\Pi$ denote the set of all prime numbers and define the
function $\sigma:\Pi \to \Pi$ by

\begin{equation*}
	\sigma(p)=\begin{cases}
		2p+1 & \text{if } 2p+1\in \Pi \\
		 p & \text{ otherwise.}
	\end{cases}
\end{equation*} \vspace{2mm}

Let $\Pi'=\{p\in \Pi: (p-1)/2 \notin \Pi \}$ and  
$\Pi''=\{p\in \Pi': 2p+1 \notin \Pi \}$.
It is easy to show that the set $\Pi'$ is infinite. Indeed
if $p \in \Pi$ is congruent to 1 modulo 4 then the integer 
$(p-1)/2$ is even and so it is a prime only if $p=5$. In all 
other cases $p \in \Pi'$. 
For $p \in \Pi'$ we set $\Sigma_p=\{\sigma^i(p):i \ge 0\}$ 
and we note that $\Pi' \cap \Sigma_p=\{p\}$.
Evidently the sets $\Sigma_p$ with $p\in \Pi'$ form a 
partition of $\Pi$. We shall refer to the sets $\Sigma_p$
with $p \in \Pi'$ as {\em prime chains} 
and we say that $p$ is the {\em head} of $\Sigma_p$. 
A prime chain $\Sigma_p$ is {\em trivial} if $\Sigma_p=\{p\}$.

We are indebted to H. Kharaghani for informing us that
the prime chains $\Sigma_p$ are well known in number theory.
They are the same as {\em the complete Cunningham chains of the first type}.
It is known that they are all finite sets. 

A prime $p$ is a {\em Sophie Germain prime} if $2p+1$ is also a prime. 
It is conjectured (but still not proved) that there exist 
infinitely many Sophie Germain primes.
This is equivalent to the claim that there are infinitely many nontrivial $\Sigma_p$.
Although there are several other conjectures (and many computational results) about $\Sigma_p$,
we shall propose one more. 

Let $\pi_N$ be the number of primes $\le N$,
$\pi'_N$ the number of prime chains $\Sigma_p$ with head $p\le N$,
and $\pi''_N$ the number of these chains which are also trivial.
The table below summarizes the results of our computations.
For each value of the exponent $e=1,2,\ldots,8$, 
the number $N$ is the smallest prime larger than $10^e$.
The timings in seconds are given just for computing the data for the primes
not covered by the previous cases (smaller values of $e$). \vspace{2mm}

$ \begin{array}{rrrrrlr}
	e & N \hspace{7mm} & \pi_N \hspace{5mm} & \pi'_N \hspace{4mm} & \pi''_N \hspace{3mm} & \hspace{3mm} {\pi''_N}/{\pi'_N} & {\rm seconds} \\[1mm]
	\hline \\[-3mm]
	1 &        11 &       5 &       2 &       0 & 0            &   0.075 \\
        2 &       101 &      26 &      19 &      13 & 0.6842105263 &   0.074 \\
	3 &      1009 &     169 &     144 &     114 & 0.7916666667 &   0.072 \\
	4 &     10007 &    1230 &    1114 &     952 & 0.8545780969 &   0.135 \\
	5 &    100003 &    9593 &    8923 &    7878 & 0.8828869214 &   0.616 \\
	6 &   1000003 &   78499 &   74175 &   67135 & 0.9050893158 &   5.016 \\
	7 &  10000019 &  664580 &  633923 &  582143 & 0.9183181554 &  52.412 \\
	8 & 100000007 & 5761456 & 5531888 & 5136716 & 0.9285647143 & 569.637 
\end{array} $ \vspace{2mm}

We conjecture that almost all $\Sigma_p$ are trivial, i.e. that
$$\lim_{N\to\infty} \frac{\pi''_N}{\pi'_N}=1.$$

\section{Acknowledgements}
This research was enabled in part by support provided by SHARCNET (http:// 
\\ www.sharcnet.ca) and 
the Digital Research Alliance of Canada (alliancecan.ca).


\begin{thebibliography}{99}

\bibitem{BDKG:best2019}
C. Bright, D. {\v{Z}}. {\Dbar}okovi{\'c}, 
I. Kotsireas, V. Ganesh,
The SAT+CAS method for combinatorial search with applications 
to best matrices,
Annals of Mathematics and Artificial Intelligence 
(2019) 87:321--342.

\bibitem{Djokovic:JCMCC:1993}
D. {\v{Z}}. {\Dbar}okovi{\'c},
Good Matrices of Orders 33,35 and 127.
JCMCC 14 (1993), 145--152.

\bibitem{DGK:JCD:2014}
D. {\v{Z}}. {\Dbar}okovi{\'c}, O. Golubitsky and 
I. S. Kotsireas,
Some new orders of Hadamard and skew-Hadamard matrices.
J. Combin. Designs, 22 (2014), 270--277.

\bibitem{DK:Goethals-Seidel}
D. {\v{Z}}. {\Dbar}okovi{\'c}, I. S. Kotsireas,
Goethals-Seidel difference families with symmetric or skew 
base blocks.
Math. Comput. Sci. (2018) 12: 373--388.

\bibitem{GS-array:1970}
J. M. Goethals and J. J. Seidel,
A skew-Hadamard matrix of order 36,
J. Austral. Math. Soc. A 11 (1970), 343--344.

\bibitem{HKT}
W. H. Holzmann,  H. Kharaghani and B. Tayfeh-Rezaie, 
Williamson matrices up to order 59,
Designs, Codes and Cryptography, 46 (2008),  343--352.


\bibitem{CK-Had:2007}
R. Craigen and H. Kharaghani,
Hadamard matrices and Hadamard designs, 
in Handbook of Combinatorial Designs, 2nd ed. C. J. Colbourn, J. H. Dinitz (eds)  pp. 273--280.
Discrete Mathematics and its Applications (Boca Raton).
Chapman \& Hall/CRC, Boca Raton, FL, 2007.

\bibitem{KT:428:2005}
H. Kharaghani and B. Tayfeh-Rezaie,
A Hadamard matrix of order 428,
J. Combin. Des. 13 (6): 435-440, 2005.

\bibitem{KouStyl:2008}
C. Koukouvinos, S. Stylianou, 
On skew-Hadamard matrices,
Discrete Math. 308 (2008), 2723--2731.

\bibitem{Mpl}
Maple 2023 (X86 64 LINUX) Maplesoft, a division of Waterloo Maple Inc.,
Waterloo, Ontario.

\bibitem{SY:1992} J. Seberry, M. Yamada, Hadamard matrices,
sequences, and block designs. In Contemporary design theory,
431-560, Wiley-Intersci. Ser. Discrete Math. Optim.,
Wiley, New York, 1992.

\bibitem{SY:2020} J. Seberry, M. Yamada, Hadamard matrices,
Constructions using Number Theory and Algebra,
Wiley, New York, 2020.

\bibitem{LW:Combinatorics:1992} J. H. van Lint \& R. M. Wilson,
A course in Combinatorics, Cambridge University Press, 1992

\end{thebibliography}
\end{document}